\documentclass[a4paper,twoside,11pt, reqno]{amsart}
\newcommand{\Ueberschrift} {Cuspidal sections for varieties of arbitrary dimension}
\newcommand{\Kurztitel}{Cuspidal sections for varieties of arbitrary dimension}

\usepackage[margin=1.5in]{geometry}
\usepackage[headings]{fullpage}
\DeclareRobustCommand{\SkipTocEntry}[5]{}
\usepackage{amsmath}
\usepackage{amsthm}
\usepackage{mathtools}
\usepackage{mathrsfs}
\usepackage{amssymb}
\usepackage{nicefrac}
\usepackage{xcolor}
\definecolor{love}{RGB}{128, 15, 37}
\definecolor{slate}{RGB}{42, 54, 59}
\usepackage{stix2}
\usepackage[LGR,T1]{fontenc} 
\usepackage[utf8]{inputenc} 
\usepackage{url} 
\usepackage{hyperref}
\definecolor{darkblue}{rgb}{0,0,.85} 
\definecolor{darkred}{rgb}{0.84,0,0}
\hypersetup{
  colorlinks = true, 
  linktoc = all,     
  linkcolor = love,  
  citecolor = love,
  urlcolor=darkblue,       
  pdfencoding = unicode,
  hypertexnames = true,
  bookmarksdepth = subsubsection,
}

\usepackage[alphabetic,nobysame]{amsrefs}
\usepackage{csquotes}
\usepackage[noabbrev,nameinlink,capitalize]{cleveref} 
\newcommand{\stackstag}[1]{\href{https://stacks.math.columbia.edu/tag/#1}{Tag #1}}
\newcommand{\stacks}[1]{\cite[\stackstag{#1}]{stacks-project}}
\numberwithin{equation}{section} 

\theoremstyle{definition}

\newtheorem{definition}[equation]{Definition}

\newtheorem{remark}[equation]{Remark}

\newtheorem*{slogan*}{Slogan}

\theoremstyle{plain}
\newtheorem{conjecture}[equation]{Conjecture}

\newtheorem{lemma}[equation]{Lemma}
\newtheorem{proposition}[equation]{Proposition}
\newtheorem{theorem}[equation]{Theorem}

\newtheorem{thmABC}{Theorem}

\newtheorem{propABC}[thmABC]{Proposition}

\crefname{section}{\S\!}{\S\S\!}
\Crefname{section}{Section}{Sections}
\crefname{subsection}{\S\!}{\S\S\!}
\Crefname{subsection}{Subsection}{Subsections}

\crefformat{nul}{(#2#1#3)} 
\crefname{nul}{}{}
\Crefformat{nul}{(#2#1#3)}
\Crefname{nul}{}{}

\crefname{axiom}{Axiom}{Axioms}
\Crefname{axiom}{Axiom}{Axioms}
\crefname{convention}{Convention}{Conventions}
\Crefname{convention}{Convention}{Conventions}
\crefname{conjecture}{Conjecture}{Conjectures}
\Crefname{conjecture}{Conjecture}{Conjectures}
\crefname{construction}{Construction}{Constructions}
\Crefname{construction}{Construction}{Constructions}
\crefname{corollary}{Corollary}{Corollaries}
\Crefname{corollary}{Corollary}{Corollaries}
\crefname{counterexample}{Counterexample}{Counterexamples}
\Crefname{counterexample}{Counterexample}{Counterexamples}
\crefname{definition}{Definition}{Definitions}
\Crefname{definition}{Definition}{Definitions}
\crefname{example}{Example}{Examples}
\Crefname{example}{Example}{Examples}
\crefname{exercise}{Exercise}{Exercises}
\Crefname{exercise}{Exercise}{Exercises}
\crefname{lemma}{Lemma}{Lemmas} 
\Crefname{lemma}{Lemma}{Lemmas} 
\crefname{notation}{Notation}{Notations}
\Crefname{notation}{Notation}{Notations}
\crefname{observation}{Observation}{Observations}
\Crefname{observation}{Observation}{Observations}
\crefname{proposition}{Proposition}{Propositions}
\Crefname{proposition}{Proposition}{Propositions}
\crefname{question}{Question}{Questions}
\Crefname{question}{Question}{Questions}
\crefname{recollection}{Recollection}{Recollections}
\Crefname{recollection}{Recollection}{Recollections}
\crefname{remark}{Remark}{Remarks}
\Crefname{remark}{Remark}{Remarks}
\crefname{remarks}{Remarks}{Remarks}
\Crefname{remarks}{Remarks}{Remarks}
\crefname{subexample}{Subexample}{Subexamples}
\Crefname{subexample}{Subexample}{Subexamples}
\crefname{theorem}{Theorem}{Theorems}
\Crefname{theorem}{Theorem}{Theorems}
\crefname{variant}{Variant}{Variants}
\Crefname{variant}{Variant}{Variants}
\crefname{warning}{Warning}{Warnings}
\Crefname{warning}{Warning}{Warnings}
\usepackage[shortlabels]{enumitem}
\setlist[enumerate]{
     itemsep  = 0.2cm,
       label  = {\upshape (\arabic*)},
         ref  = \arabic*,
    leftmargin  = *
}

\newcommand{\enumref}[2]{{\upshape (\hyperref[#1.#2]{\ref*{#1}.\ref*{#1.#2}})}}
\usepackage{xspace}
\usepackage[final,factor=1100,stretch=10,shrink=10]{microtype}
\SetProtrusion{encoding={*},family={bch},series={*},size={6,7}}
              {1={ ,750},2={ ,500},3={ ,500},4={ ,500},5={ ,500},
               6={ ,500},7={ ,600},8={ ,500},9={ ,500},0={ ,500}}

\SetExtraKerning[unit=space]
    {encoding={*}, family={bch}, series={*}, size={footnotesize,small,normalsize}}
    {\textendash={400,400},   
     "28={ ,150},     
     "29={150, },     
     \textquotedblleft={ ,150},   
     \textquotedblright={150, }}  

\SetTracking{encoding={*}, shape=sc}{40}
\usepackage{tikz}
\usepackage{tikz-cd}
\usepackage{amscd}

\usepackage[all]{xy} \CompileMatrices
\usepackage{extarrows}
\tikzcdset{arrow style=tikz, diagrams={>=stealth}} 

\newcommand{\colonequals}{\coloneq}

\renewcommand{\setminus}{\smallsetminus}

\DeclareMathOperator{\Hom}{Hom}

\DeclareMathOperator{\im}{im}

\DeclareMathOperator{\pr}{pr}
\DeclareMathOperator{\Gal}{Gal}

\DeclareMathOperator{\Spec}{Spec}

\newcommand{\spec}{\Spec}
\newcommand{\homo}{\Hom}
\newcommand{\cl}[1]{\overline{#1}}
\newcommand{\Oc}{\mathcal{O}}
\newcommand{\Zb}{\bZ}

\newcommand{\sections}[2]{\mathscr{S}_{\pi_1(#1/#2)}}
\newcommand{\seck}[1]{\sections{#1}{k}}

\newenvironment{defn}
  {\begin{definition}}{\end{definition}}
\newenvironment{rem}
  {\begin{remark}}{\end{remark}}
\newenvironment{conj}
  {\begin{conjecture}}{\end{conjecture}}
\newenvironment{prop}
  {\begin{proposition}}{\end{proposition}}
\newenvironment{thm}
  {\begin{theorem}}{\end{theorem}}

\newcommand{\bZ}{{\mathbb Z}}

\newcommand{\cC}{{\mathscr C}}

\renewcommand{\epsilon}{\varepsilon}
\renewcommand{\emptyset}{\varnothing}

\begin{document}

\title[\Kurztitel]{\Ueberschrift} 

\author{Magnus Carlson}
\address{Magnus Carlson
\newline \indent
Department of Electrical Engineering, Mathematics and Science
\newline \indent
University of Gävle, Sweden}
\email{\tt magnus.carlsson@hig.se}

\author{David Corwin}
\address{David Corwin
\newline \indent
Department of Mathematics, Ben-Gurion University of the Negev
\newline \indent
P.O.B. 653, Be'er Sheva 8410501, Israel}
\email{\tt corwind@bgu.ac.il}

\date{\normalsize \today}

\begin{abstract} 
We define cuspidal sections for varieties of arbitrary dimension over fields of characteristic zero, extending the definition for curves by Grothendieck. Our notion of cuspidal sections is defined relative to a normal crossings compactification. We prove that the closure of the resulting set is independent of the choice of a compactification. For varieties iteratively constructed from fibrations, we compare our definition to an inductive notion of cuspidal sections, and discuss connections to the section conjecture of Grothendieck.  \end{abstract}

\maketitle

\section{Introduction}
The \emph{section conjecture} of Grothendieck \cite{Gro83} posits that rational points of proper hyperbolic curves $U$ over number fields $k$ are in canonical bijection with (conjugacy classes) of sections of the natural map $\pi_1(U) \to G_k$. If $U$ is hyperbolic but not proper, there are often more sections of $\pi_1(U) \to G_k$ than there are rational points. To account for this, Grothendieck introduced \emph{cuspidal sections} which are associated to the rational points at infinity. The extension of Grothendieck's section conjecture to (not necessarily proper) hyperbolic curves is then that any section of $\pi_1(U) \to G_k$ either comes from a rational point, or is cuspidal. 

In this paper, we first define, given a smooth geometrically connected variety $U$ over $k$ and a normal crossings compactification $X$ of $U$, an associated set $C_{U,X}$ of cuspidal sections (Definition \ref{defn:higher_dim_cusps}), generalizing Grothendieck's definition for curves. We prove that while $C_{U,X}$ depends on the compactification, its closure $C_U$ does not. When $U$ is a hyperbolic polycurve, a variety built out of hyperbolic curves via fibrations, we relate the set $C_U$ with a related notion of cuspidal sections on $U$, built inductively from the fibrations defining $U$. We now state our results more precisely.

Let $k$ be a field of characteristic zero with a given algebraic closure $k^s$ and let $U/k$ be a smooth geometrically connected variety. We set $G_k = \Gal(k^s/k)$ and write $\seck{U}$ for the set of sections of the surjective map $\pi_1(U) \to G_k$ modulo conjugacy by the geometric fundamental group $\pi_1(U^s)$, where $U^s = U \times_k k^s$. Given a normal crossings compactification $U \subset X$, we obtain a set $C_{U,X} \subset \seck{U}$ of cuspidal sections: for each $y \in (X \setminus U)(k)$, we let $X_{(y)}$ be the Henselization of $X$ at $y$ and $U_{(y)} = X_{(y)} \times_X U$. For each $y$, there is a map $\seck{U_{(y)}} \to \seck{U}$; $C_{U,X}$ is the union of these images for $y \in (X \setminus U)(k)$. Recall that the set of sections $\seck{U}$ has a natural topology (see \cite[\S 4.2]{sxbook13}).
\begin{thmABC}[{Theorem \ref{thm:cuspidals_independent_of_compactification}}] 
Let $k$ be a field of characteristic zero, and let $U/k$ be a smooth variety and let $U \hookrightarrow X$ be a normal crossings compactification. Then the closure $C_U$ of the subset $C_{U,X}$ does not depend on the normal crossings compactification $X$.
\end{thmABC}
A hyperbolic polycurve $X/k$ is a geometrically connected variety $X$ admitting a sequence of parameterizing morphisms 
\[ X = X_n \xrightarrow{f_n} X_{n-1} \cdots \xrightarrow{f_2} X_1 \xrightarrow{f_1} \Spec k \] such that each $X_i \to X_{i-1}$ is a relative hyperbolic curve (see \cite[Definition 2.1]{Hoshi14}). Given a hyperbolic polycurve $X/k$ and a sequence of parameterizing morphisms above, one inductively constructs a cuspidal set $C_{X,f} \subset \seck{X}$: at each stage, a section is cuspidal if its projection to the base is cuspidal or if it comes from a cuspidal section of a fiber over a rational point. Since any hyperbolic curve has a canonical cuspidal set, this gives a well-defined set of cuspidal sections $C_{X,f}$, a priori depending on the sequence of parameterizing morphisms. We say that a hyperbolic curve $X/k$ \emph{satisfies injectivity in the section conjecture} if the non-abelian Kummer map $\kappa_{X/k} \colon X(k) \to \seck{X}$ is injective with image disjoint from the cuspidal sections defined by Grothendieck, and the same holds for every geometrically connected finite étale cover $X' \to X$.
\begin{propABC}[{Proposition \ref{prop:equivsection}}]
Let $k$ be a field of characteristic zero, and assume that every hyperbolic curve over $k$ satisfies injectivity in the section conjecture. Then the cuspidal set $C_{X,f}$ is independent of the sequence of parameterizing morphisms for every hyperbolic polycurve $X$ over $k$ if and only if the section conjecture holds for hyperbolic curves over $k$.
\end{propABC}
We further investigate the relationship between the inductively constructed $C_{X,f}$ and the set of cuspidal sections $C_X$ obtained as the closure of the set of cuspidal sections associated to a normal crossings compactification.  
\begin{propABC}[{Proposition \ref{prop:fibinclusion}}]
Let $k$ be a field of characteristic zero, and let $X$ be a hyperbolic polycurve over $k$. Then $C_{X} \subset \overline{C_{X,f}}$.
\end{propABC}
\subsection{Notation and Conventions}\label{sec:noteconv}
Whenever we speak of a field $k$, we implicitly fix a separable closure $k^s$ throughout, which produces a geometric point $\spec(k^s) \to \spec(k)$ of $k$. A geometric point $\overline{x}$ of a scheme $X$ is a map $\Spec \Omega \rightarrow X$, where $\Omega$ is a separably closed field.  A \emph{variety} over $k$ is a finite type $k$-scheme which is reduced, separated, and geometrically connected. A \emph{curve} is a smooth (geometrically connected) variety of dimension one. For a scheme $X$ over a field $k$ we denote by $X^s$ the base change of $X$ to $k^s$. If $X$ is a scheme, $\Oc(X)$ denotes the global sections $\Oc_X(X)$. For a scheme $X$ and a point $y \in X$, we denote by $X_{(y)}$ the spectrum of the Henselization of $X$ in $y$. Similarly, given a geometric point $\overline{y}$ of $X$, $X_{(\overline{y})}$ is the \emph{strict} Henselization of $X$ in $\overline{y}$.  

When the choice of a basepoint is irrelevant, we write $\pi_1(X)$ for the étale fundamental group of $X$ fixed at some unspecified basepoint. A homomorphism $\varphi \colon \pi_1(X) \to \pi_1(Y)$ is then understood to respect the choices of basepoints. Given a map of schemes $f \colon X \to Y$, we write $f_* \colon \pi_1(X) \to \pi_1(Y)$ for the induced map of étale fundamental groups.
\subsection{Acknowledgements}
The second author would like to thank Jakob Stix for suggesting that this be separate from \cite{OnPoonenCE} and Tomer Schlank for sharing the idea of cuspidal sections for certain higher-dimensional varieties (specifically, for good neighborhoods) as part of the work on \cite{OnPoonenCE}, which resulted in this paper. He would like to thank Olivier Wittenberg, for help with the overall approach to proving independence of compactification and hosting at ENS, Kenji Matsuki and Dan Abramovich for answering numerous questions about resolution of singularities, and Martin Olsson for numerous discussions about logarithmic fundamental groups. Finally, he would like to thank his doctoral advisor Bjorn Poonen for an enormous amount of help, including multiple careful readings of this and other related papers, as well as numerous discussions about the ideas.

\section{Background on the Section Conjecture}\label{sec:grosec_conj}

We review the section conjecture of Grothendieck.

Fix a separable closure $k^s$ of $k$, which produces a geometric point $\spec(k^s) \to \spec(k)$ of $k$. We identify $G_k  \colonequals  \mathrm{Gal}(k^s/k)$ with the \'etale fundamental group of $\spec{k}$ based at this point. 

Given a qcqs geometrically connected scheme $X$ over $k$, and a geometric point $\overline{x}$ of $X^s$, there is the \emph{fundamental exact sequence} (see \stacks{0BTX})
\[ 1 \to \pi_1(X^s,\cl{x}) \to \pi_1(X,\overline{x}) \to G_k \to 1. \] 

\begin{defn}\label{defn:seck}
Let $X$ be a geometrically connected qcqs scheme over $k$. We define the set of sections $\mathscr{S}_{\pi_1(X/k),\cl{x}}$ as the set of sections of the surjection $\pi_1(X,\cl{x}) \to G_k$ modulo the action of $\pi_1(X^s,\cl{x})$ by conjugation.\end{defn}
As is well-known \cite[Definition 23]{sxbook13}, the sets  $\mathscr{S}_{\pi_1(X/k),\cl{x}}$ for different choices of $\cl{x}$ are in canonical bijection, so we may write $\seck{X}$ without ambiguity.







A map $ f \colon X \to Y $  of geometrically connected qcqs $k$-schemes induces a map $f_* \colon \pi_1(X,\overline{x}) \to \pi_1(Y,\overline{y})$, well-defined up to conjugation by $\pi_1(Y^s,\overline{y})$. We also denote the induced map $\seck{X} \to \seck{Y}$ by $f_*$, by abuse of notation. The \emph{profinite Kummer map} is the map $\kappa_{X/k} \colon X(k) \to \seck{X}$ induced by this functoriality of $\pi_1$ up to conjugation.

\begin{conj}[Grothendieck \cite{Gro83}]\label{c:gro1} Let $k$ be a number field
and $X/k$  a smooth proper curve of genus at least $2$. Then the Kummer map $X(k) \to \seck{X}$ is a bijection.
\end{conj}
\begin{rem}
The above conjecture can be extended to other base fields. OpenAI has recently posted a claimed proof of the analogue of Conjecture \ref{c:gro1} when $k$ is a $p$-adic local field \cite{OpenAI2026Section}. 
\end{rem}

Grothendieck, in his letter to Faltings, also extended this conjecture to the non-proper case \cite{Gro83}. This extended conjecture requires the notion of cuspidal sections.


\subsection{Cuspidal Sections for Curves}
\begin{defn}[{see \cite{Gro83}}]\label{defn:cusp_curve}
Let $k$ be a field of characteristic zero, and let $U$ be a smooth curve over $k$ with smooth compactification $X$, and set $Y = X - U$. Given $y \in Y(k)$, we define $C_y$, \emph{the packet of cuspidal sections at $y$}, as the image of the map $\seck{U_{(y)}} \to \seck{U}$. 
The set $C_U$ of all \emph{cuspidal sections of the curve $U$} is defined to be
\[  \bigcup_{y \in Y(k)} \mathrm{Im}\left(\seck{U_{(y)}} \to \seck{U}\right) \subseteq \seck{U}.\]
\end{defn}


A curve $U$ over $k$ is \emph{hyperbolic} if its Euler characteristic over $k^s$ is negative, or, equivalently, $\pi_1(U^s)$ is non-abelian. The following extends Conjecture \ref{c:gro1} to non-proper hyperbolic curves:





\begin{conj}[Grothendieck \cite{Gro83}]\label{c:gro2} Let $k$ be a number field
and suppose that $U$ is hyperbolic. Then every section $s \in \seck{U}$ is either in the image of $\kappa_{U/k}$, or is cuspidal.
\end{conj}

\section{Cuspidal Sections in Arbitrary Dimension}\label{sec:cusp_arb_dim}
\subsection{Preliminaries on Schemes of Nearby Points}
We define the scheme of nearby points; a similar definition can be found in \cite[Chapter 18]{sxbook13}.
\begin{defn}\label{defn:nearby_points}
Let $U \subseteq X$ be an open immersion of schemes. For $y \in X$ (resp. for a geometric point $\overline{y} \to X$), the \emph{scheme of nearby points to $U$ at $y$} (resp. the scheme of \emph{strict nearby points to $U$ at $\overline{y}$})  is
\[ U_{(y)} := X_{(y)} \times_X U \qquad (\text{resp. } U_{(\overline{y})} := X_{(\overline{y})} \times_X U)  \]  where $X_{(y)}$ (resp. $X_{(\overline{y})}$) is the Henselization of $X$ (resp. strict Henselization of $X$) at $y$.
\end{defn}
Suppose that $U \subseteq X$ is a dense open immersion of $k$-varieties  and $y \in X(k)$. Then, if $X$ is geometrically unibranch at $y$, $U_{(y)}$ is geometrically irreducible over $k$ \stacks{06DM}, hence geometrically connected.
\begin{defn}\label{defn:cuspidal_packet}
Let $U \hookrightarrow X$ be an open immersion of smooth $k$-varieties. Let $Y = X - U$ and $y \in Y(k)$. The \emph{packet of cuspidal sections at $y$}, denoted $C_y$, is the image of the natural map \[ \seck{U_{(y)}} \to \seck{U},\] associated to the canonical map $U_{(y)} \to U$.
\end{defn}

\begin{lemma}
\label{lemma:excision}
Let $U \to X$ be an open immersion of $k$-schemes, $f \colon X' \to X$ a map of $k$-schemes and set  $U' = f^{-1}(U)$. Let $y \in X'(k)$ (resp. $\overline{y} \in X'(k^s)$). Then $f$ induces a natural map\footnote{with respect to the maps to $U'$ and $U$} \[ U'_{(y)} \to U_{(f(y))} \]  (resp. $U'_{(\overline{y})} \to U_{(f(\overline{y}))}$), which is an isomorphism if $f$ is \'etale.\end{lemma}

\begin{proof}
By functoriality of henselization, we obtain a map $X'_{(y)} \to X_{(f(y))}$. As $U' = X' \times_X U$, we obtain by base-change a functorial map $U'_{(y)} \to U_{(f(y))}$. 

When $f$ is étale, the morphism $X'_{(y)} \to X_{(f(y))}$ is an isomorphism since their residue fields are isomorphic. Hence, the base-change $U'_{(y)} \to U_{(f(y))}$ is an isomorphism as well.

The claim for schemes of strict nearby points is proven analogously.
\end{proof}

\subsection{Cuspidal Sections Associated to a Compactification}
We first calculate the étale fundamental group of the scheme of nearby points $U_{(y)}$ for a normal crossings compactification $j \colon U \to X$.
\begin{prop} \label{prop:fundamentalofU}
Let $k$ be a field of characteristic zero, and let $U$ be a smooth $k$-variety and $j \colon U \to X$ be a normal crossings compactification.  Set $Y = X-U$ and let $y \in Y(k)$, and write $\overline{y} \to Y$ for the geometric point induced by $y$ and our fixed separable closure $k^s/k$. We then have a short exact sequence
\[
1 \to \pi_1(U_{(\overline{y})}) \cong \homo(\Oc(U_{(\overline{y})})^\times/\Oc(X_{(\overline{y})})^\times,\widehat{\Zb}(1)) \to \pi_1(U_{(y)}) \to G_k \to 1,
\] 
and the action of $G_k$ is obtained from the natural actions on $\Oc(U_{(\overline{y})})^\times/\Oc(X_{(\overline{y})})^\times$ and $\widehat{\Zb}(1)$.
This exact sequence is functorial in the pair $(X,U)$.
\end{prop}

\begin{proof}
Let $(X_{(y)},M)$ be the log scheme associated to the divisor defined by $Y$ in $X_{(y)}$ and $(X_{(\overline{y})},M^s)$ the pullback of the log structure along $X_{(\overline{y})} \to X_{(y)}$. Thus, the log structure on $X_{(y)}$ is obtained by pullback from the sheaf of monoids $j_*\Oc_U^{\times} \cap \Oc_X$ on the small étale site of $X$. Write $\overline{M} \colonequals M/\Oc_{X_{(y)}}^\times$ for the associated characteristic sheaf of monoids on $X_{(y)}$. We let $s$ be the log-point $(y,M_y)$, obtained from the log structure of $X_{(y)}$ by pullback along  $y \colon \Spec k \to X_{(y)}$, and denote by $\overline{s} = s \otimes_k k^s$ the log-point obtained via base change.

By \cite[Examples 4.7(c), (d)]{IllusieOverview} (see also \cite[Remark 10.3]{KatoLogII}), we have isomorphisms  \[ \pi_1(U_{(y)}) \cong \pi_1^{\log}(X_{(y)}) \cong \pi_1^{\log}(s), \] and moreover, these isomorphisms are compatible with the maps to $G_k$. In the same way, we obtain an isomorphism $\pi_1(U_{(\overline{y})}) \cong \pi_1^{\log}(\overline{s})$, compatible with the maps $\pi_1(U_{(\overline{y})}) \to \pi_1(U_{(y)})$ and $\pi_1^{\log}(\overline{s}) \to \pi_1^{\log}(s)$.

Since $Y$ is normal crossings, we obtain, by taking orders along local (geometric) boundary components, that the group completion of the characteristic monoid $\overline{M}_{\overline{y}}$ is isomorphic to $\Oc(U_{(\overline{y})})^\times/\Oc(X_{(\overline{y})})^\times$ as a $G_k$-module.

By \cite[Example 4.7(a)]{IllusieOverview}, we then have an associated short exact sequence
\[
1 \to \homo(\Oc(U_{(\overline{y})})^\times/\Oc(X_{(\overline{y})})^\times,\widehat{\Zb}(1)) \to \pi_1^{\log}(s) \to G_k \to 1,
\]
where the kernel is identified with $\pi_1^{\log}(\overline{s})$. These comparisons, together with their naturality in the pair $(X,U)$, imply our proposition.
\end{proof}

\begin{rem}
A similar sequence appears in \cite[(18.3)]{sxbook13} but contains a typo: its use of nearby points in the kernel should be replaced by strict nearby points.
\end{rem}



\begin{defn}\label{defn:higher_dim_cusps}

Let $U$ be a smooth geometrically connected variety over a field $k$ of characteristic $0$, let $U \hookrightarrow X$ be a normal crossings compactification and set $Y=X \setminus U$.


We define the set $C_{U,X}$ of all \emph{cuspidal sections of the variety $U$ relative to $X$} as
$$
C_{U,X} := \bigcup_{y \in Y(k)} \mathrm{Im}\left(\seck{U_{(y)}} \to \seck{U}\right) \subseteq \seck{U}.
$$
\end{defn}
\begin{thm}\label{thm:cuspidals_independent_of_compactification}
Let $k$ be a field of characteristic zero, and let $U/k$ be a smooth variety and let $U \hookrightarrow X$ be a normal crossings compactification. Then the closure $C_U$ of the subset $C_{U,X}$ does not depend on the normal crossings compactification $X$.
\end{thm}
\begin{remark}
When $U/k$ is a curve, any normal compactification $U \hookrightarrow X$ is \emph{the} unique smooth compactification of $U$. In that case $C_{U,X}$ is automatically closed, thus calling this closure $C_U$ is consistent with Definition \ref{defn:cusp_curve}. To see this, we can without loss of generality assume that $\pi_1(U^s)$ is non-trivial, thus the Euler characteristic $\chi(U^s) \leq 0$. Let $\widetilde{U} \to U$ be the universal pro-étale cover of $U$, and $\widetilde{X}$ be the normalization of $X$ in $\widetilde{U}$, and $\widetilde{Y}$ the inverse image of $Y$ in $\tilde{X}$, as in \cite[\S 3.1]{sx08}. We let $\widetilde{\operatorname{Cusps}}(U) = \Hom_{k^s}(\Spec k^s,\widetilde{Y})$ be the set of prolongations of the cusps of $U$ to $\widetilde{U}$: this is a profinite set and has a continuous $\pi_1(U)$ action (cf. \cite[Definition 10]{sx08}). It follows from \cite[\S 3.3 and Lemma 14]{sx08} that a section $s \in \seck{U}$ is cuspidal if and only if some (equivalently, any) representative fixes an element of $\widetilde{\operatorname{Cusps}}(U)$. Since $\widetilde{\operatorname{Cusps}}(U)$ is a compact set, the set of closed subgroups of $\pi_1(U)$ that fix an element of $\widetilde{\operatorname{Cusps}}(U)$ is a closed subset of $\operatorname{Sub}(\pi_1(U))$, the set of closed subgroups of $\pi_1(U)$, where the topology is that of \cite[\S 4.1]{sxbook13}. As this subset is invariant under conjugation by $\pi_1(U^s)$, the description of the topology on $\seck{U}$ in \cite[\S 4.1]{sxbook13} shows that $C_{U,X}$ is closed. \end{remark}
The rest of this section will be concerned with proving Theorem \ref{thm:cuspidals_independent_of_compactification}. We first state and prove the following well-known lemma.
\begin{lemma}\label{lemma:birat_factorization}
Let $U$ be a smooth geometrically connected variety over a field $k$ of characteristic zero, and let $X,X'$ be two smooth compactifications of $U$. Then there is an snc compactification $X''$ along with maps $X'' \to X'$ and $X'' \to X$ that are isomorphisms over $U$.
\end{lemma}

\begin{proof}
Let $\Delta$ be the diagonal of $U$ in $U \times U$, and let $\cl{\Delta}$ be its closure in $X' \times X$. Then $\cl{\Delta}$ is a proper but possibly singular variety, and it has morphisms to both $X'$ and $X$. Now $\Delta$ is closed in $U \times X$ and $X' \times U$, respectively, because it is the graph of a morphism in each case, so the intersection of $\cl{\Delta}$ with $U \times X$ and $X' \times U$, respectively, is just $\Delta$. Therefore, the projections of $\cl{\Delta}$ to $X'$ and $X$, respectively, are both isomorphisms over $U$. By resolution of singularities in characteristic $0$, there is a smooth and proper variety  $X''$ along with a proper birational map $X'' \to \cl{\Delta}$ that is an isomorphism over $\Delta \subseteq \cl{\Delta}$, such that  $X'' \setminus \Delta$ is a simple normal crossings divisor. 
\end{proof}

\begin{lemma}\label{lemma:cusps_morphism_of_compactifications}
Let $U$ be a smooth variety over a field $k$ of characteristic zero and let $X$ and $X'$ be two smooth normal crossings compactifications of $U$. If $h \colon X' \to X$ is a morphism which is an isomorphism over $U$, then  $C_{U,X'} \subseteq C_{U,X}$.
\end{lemma}

\begin{proof}
Let $y \in X'(k) \setminus U(k)$. 




By Lemma \ref{lemma:excision}, we have the following diagram,  because $X' \to X$ is an isomorphism over $U$:

$$
\begin{CD}
U_{(y)} @>>> U\\
@VVV @VVV\\
U_{(h(y))} @>>> U
\end{CD}
$$

The right vertical arrow is the identity, so by taking sections, we obtain a diagram
$$
\begin{CD}
\seck{U_{(y)}} @>>> \seck{U}\\
@VVV @VVV\\
\seck{U_{(h(y))}} @>>> \seck{U} 
\end{CD}
$$

It follows that the image of $\seck{U_{(y)}}$ in $\seck{U}$ is in the image of $\seck{U_{(h(y))}}$, hence in $C_{U,X}$. Thus $$C_{U,X'} = \bigcup_{y \in X'(k) \setminus U(k)} \mathrm{Im}\left(\seck{U_{(y)}} \to \seck{U}\right) \subseteq C_{U,X}.$$
\end{proof}
\begin{lemma}\label{lemma:cusps_propermod}
Let $U$ be a smooth variety over a field $k$ of characteristic zero, and let $X$ and $X'$ be smooth normal crossings compactifications of $U$. Let $h \colon X' \to X$ be a proper morphism that is an isomorphism over $U$. Then we have an inclusion $\overline{C_{U,X}} \subseteq \overline{C_{U,X'}}$.
\end{lemma}
\begin{proof}
Let $y \in X(k) \setminus U(k)$. We will show that 
\[ \bigcup_{y'\in h^{-1}(y)(k)}\im \left (\seck{U_{{(y')}}} \to \seck{U_{(y)}} \right) \] is dense in $\seck{U_{(y)}}$; this will imply our lemma.

Let $(X,M)$ be the log scheme associated to the normal crossings compactification $j \colon U \hookrightarrow X$; we let $\overline{M}$ be the associated monoid sheaf. Thus, we have that $M = j_*\mathcal{O}^{\times}_U \cap \mathcal{O}_X$ on the small étale site of $X$, and $\overline{M} = M/\mathcal{O}^{\times}_X$. If $\overline{M}_{\overline{y}}$ is the stalk of $\overline{M}$ at the geometric point $\overline{y}$, we let $\mathbb{T}$ be the torus over $k$ whose geometric character group is the group completion of the monoid $\overline{M}_{\overline{y}}$, with its natural $G_k$-action. Note that since our compactification is normal crossings, \cite[\S 7.3(b)]{IllusieOverview} implies that $\overline{M}_{\overline{y}}^{\operatorname{grp}}$ is the free abelian group on the set of geometric local branches at $\overline{y}$, and these branches are permuted by $G_k$.

For each $G_k$-orbit of these local branches, choose a representative $\overline{D}_i$ and let $L_i$ be its field of definition. Note that if we set $ A = \prod_i L_i$, then $\mathbb{T} = \operatorname{Res}_{A/k} \mathbb{G}_m$. Further, for each $i$, the branch divisor $\overline{D}_i$ descends to a divisor $D_i \subset X_{(y)} \times_k L_i$. Then the $D_i$'s are clearly principal; we let $x_i \in \Oc(X_{(y)}) \otimes_k L_i$ be a generator of $D_i$ and set $x = (x_i)_i \in \Oc(X_{(y)}) \otimes A$. 

We note that $x$ is a unit in $\Oc(U_{(y)}) \otimes A$, hence $x$ defines a morphism $x \colon U_{(y)} \to \mathbb{T}$. A quick computation using Proposition \ref{prop:fundamentalofU} and Kummer theory shows that the induced morphism 
\[ x_* \colon \pi_1(U_{(y)}) \to \pi_1(\mathbb{T}) \]
is an isomorphism over $G_k$. This implies that \[ \bigcup_{y'\in h^{-1}(y)(k)}\im \left (\seck{U_{{(y')}}} \to \seck{U_{(y)}} \right) \] is dense in $\seck{U_{(y)}}$ if and only if \[ \bigcup_{y'\in h^{-1}(y)(k)}\im \left (\seck{U_{{(y')}}} \to \seck{\mathbb{T}} \right) \] is dense in $\seck{\mathbb{T}}$. In addition, quasi-triviality of $\mathbb{T}$ implies by Hilbert's Theorem 90 that $\kappa_{\mathbb{T}/k}$ is surjective modulo every positive integer, hence that $\kappa_{\mathbb{T}/k}(\mathbb{T}(k))$ is dense in $\seck{\mathbb{T}}$; it therefore suffices to show that \[ \bigcup_{y'\in h^{-1}(y)(k)} x_*(\seck{U_{{(y')}}}) \] contains a translate of $\kappa_{\mathbb{T}/k}(\mathbb{T}(k))$.

To see this, let $e_1, \ldots, e_m$ be a basis of $A$ over $k$, and write \[ x = \sum_{i=1}^m z_i \otimes e_i, \qquad z_i \in \mathfrak{m}_y \Oc(X_{(y)}), \] where $\mathfrak{m}_y$ is the maximal ideal of $\Oc(X_{(y)})$. The $z_1, \ldots z_m$ form part of a  regular system of parameters for $\Oc(X_{(y)})$, so complete them to a system of parameters $z_1, \ldots, z_n$, and write $\Oc(X_{(y)}) \cong k \langle  z_1, \ldots, z_n \rangle,$ where $k\langle z_1, \ldots z_n \rangle$ is the ring of algebraic power series with coefficients in $k$. For each $a \in A^*$ write $a = \sum a_i e_i$ with $a_i \in k$, and define the map \[ t_a \colon \Spec k \langle t \rangle \to \Spec k \langle z_1, \ldots , z_n \rangle \cong X_{(y)}, \qquad t \to (a_1 t,a_2t, \ldots , a_mt,0, \ldots, 0). \] 

Let $\eta$ be the generic point of $\Spec k\langle t \rangle$ and note that the restriction of $t_a$ to $\eta$ factors through $U_{(y)}$; we denote by $\gamma_a \colon \eta \to U_{(y)}$ the restriction. We now choose an arbitrary section\footnote{exists because $\mathbb{G}_m(k) \neq \emptyset$} \[ s \colon G_k \rightarrow \pi_1(\eta), \] and consider the composite $x_* \circ (\gamma_a)_* \circ s \in \seck{\mathbb{T}}$. Since $\pi_1(\mathbb{T})$ is a group object in the category of profinite groups over $G_k$ \cite[\S 13.1]{sxbook13}, we obtain that $\seck{\mathbb{T}}$ has the structure of a group, with the identity section as the identity element. We claim that \[ x_* \circ (\gamma_a)_* \circ s = \kappa_{\mathbb{T}/k}(a)+ \alpha \in \seck{\mathbb{T}}, \]  where $\alpha = x_* \circ (\gamma_1)_* \circ s.$ Here $\gamma_1$ is the restriction of $t_1$ to $\eta$. To see this equality, note that the morphism of schemes \[ x \circ \gamma_a \colon \eta  \rightarrow \mathbb{T}\] equals $r_a \circ (x \circ \gamma_1)$, where $r_a \colon \mathbb{T} \to \mathbb{T}$ is translation by $a$. Functoriality then shows that we have the equality $x_* \circ (\gamma_a)_* \circ s = \kappa_{\mathbb{T}/k}(a)+ \alpha$. Thus, as $a$ ranges over $A^*$, the set of sections of the form $x_* \circ (\gamma_a)_* \circ s$ contains a translate of $\kappa_{\mathbb{T}/k}(\mathbb{T}(k))$.

Now, since $X' \to X$ is proper and an isomorphism over $U$, and $\eta$ lands in $U$ when we apply $t_a$, the valuative criterion shows that $t_a$ lifts to a map 
\[ t'_a \colon \Spec k \langle t \rangle  \to X'. \] 
If $y'$ denotes the image of the closed point of $\Spec k \langle t \rangle$ under $t'_a$, then $t'_a$ factors through $X'_{(y')}$. Further, letting $U' \colonequals h^{-1}(U) \cong U$, the restriction of $t'_a$ to $\eta$ factors  through $U'_{(y')}$, and we denote this restriction by $\gamma'_a$. The section $(\gamma'_a)_* \circ s \in \seck{U'_{(y')}}$ then lifts $(\gamma_a)_* \circ s$. As $a$ ranges through $A^*$, we see that the set of sections of the form $x_* \circ h_* \circ  (\gamma'_a)_* \circ s$ contains a translate of $\kappa_{\mathbb{T}/k}(\mathbb{T}(k))$. As previously noted, this implies that \[ \bigcup_{y'\in h^{-1}(y)(k)}\im \left (\seck{U_{(y')}} \to \seck{U_{(y)}} \right) \] is dense in $\seck{U_{(y)}}$.
\end{proof}
\begin{proof}[Proof of Theorem \ref{thm:cuspidals_independent_of_compactification}]
Let $U \hookrightarrow X'$ be a second compactification. By Lemma \ref{lemma:birat_factorization}, we can assume that there is a map $f \colon X' \to X$ that is an isomorphism over $U$. Then Lemmas \ref{lemma:cusps_morphism_of_compactifications} and \ref{lemma:cusps_propermod} imply that $\overline{C_{U,X}} = \overline{C_{U,X'}}$. 

\end{proof}

\section{Cuspidal Sections Induced by Fibrations}\label{sec:cuspidal_sections}
In Section \ref{sec:cusp_arb_dim}, we gave an intrinsic definition of cuspidal sections for any smooth variety $X$ over a field of characteristic zero. When $X$ is iteratively constructed from fibrations, we compare this intrinsic definition of cuspidal sections with an inductive definition of cuspidality. To aid our arguments, it is convenient to define the abstract notion of \emph{cuspidal data}. From now on, $k$ is assumed to be a field of characteristic zero.
\subsection{Cuspidal Data}
\begin{defn}\label{defn:cuspdatum}
A \emph{cuspidal datum} $(X,\cC_X)$ is a variety $X/k$ and a subset $\cC_X \subseteq \seck{X}$, called the set of \emph{cuspidal sections}.
\end{defn}
\begin{rem} \label{rem:cuspidalunderetale}
Let $X/k$ be a variety and let $f \colon X' \to X$ be a finite étale cover of varieties over $k$. Then a cuspidal datum $(X,\cC_X)$ naturally induces a cuspidal datum $(X',\cC_{X'})$ by letting \[ \cC_{X'} = f_*^{-1}(\cC_X) \subset \seck{X'}, \] where $f_* \colon \seck{X'} \to \seck{X}$ is the induced map. 
\end{rem}
\begin{defn}\label{defn:surj_sec_conj}
A cuspidal datum $(X,\cC_X)$ satisfies \emph{surjectivity in the section conjecture} if 
\[ \seck{X} = \kappa_{X/k}(X(k)) \cup \cC_X, \] that is, every section $s$ either lies in $\cC_X$ or comes from a rational point of $X$. \end{defn}
\begin{rem} \label{rem:surjectivityunderetale}
Assume that $(X,\cC_X)$ is a cuspidal datum which satisfies surjectivity in the section conjecture. Then for every finite étale cover $f \colon X' \to X$ of varieties over $k$, the induced datum $(X',\cC_{X'})$ also satisfies surjectivity in the section conjecture. Indeed, let $s \in \seck{X'}$. If $f_*(s) \in \cC_X$, then trivially $s \in \cC_{X'}$. If not, $f_*(s) \in \kappa_{X/k}(X(k)),$ and then \cite[Corollary 34]{sxbook13} shows that $s \in \kappa_{X'/k}(X'(k))$.
\end{rem}
For a cuspidal datum $(X,\cC_X)$, define the extended non-abelian Kummer map 
\[ \kappa_{X/k} \colon X(k) \amalg \cC_X \to \seck{X} \] as the ordinary non-abelian Kummer map on $X(k)$, and as the inclusion on $\cC_X$.
\begin{defn}\label{defn:inj_sec_conj}A cuspidal datum $(X,\cC_X)$ is said to satisfy \emph{injectivity in the section conjecture} if for any geometrically connected finite étale cover $X' \to X$, the extended non-abelian Kummer map
\[ \kappa_{X'/k} \colon X'(k) \amalg \cC_{X'} \to \seck{X'} \]  is injective.
\end{defn}
\begin{rem} \label{remark:cuspidalequiv}
A cuspidal datum $(X,\cC_X)$ satisfies injectivity in the section conjecture if and only if the extended Kummer map 
\[ \kappa_{X/k} \colon X(k) \amalg \cC_X \to \seck{X} \] is injective, and for any $a \in X(k)$, the centralizer $Z_{\pi_1(X^s)}(\kappa_{X/k}(a))$ is trivial\footnote{The centralizer of a conjugacy class of a section is the centralizer of the image of any representative, thus well-defined up to conjugation.}. 

If $(X,\cC_X)$ satisfies injectivity in the section conjecture, the first paragraph of \cite[\S 9.4]{sxbook13} shows that for any $a \in X(k)$, $Z_{\pi_1(X^s)}(\kappa_{X/k}(a))$ is trivial. Conversely, assume that the map \[X(k) \amalg \cC_X \to \seck{X}\] is injective and that $Z_{\pi_1(X^s)}(\kappa_{X/k}(a))$ is trivial for any $a \in X(k)$. The proof of \cite[Proposition 33(2)]{sxbook13}, together with \cite[Corollary 34]{sxbook13}, implies that for any geometrically connected finite étale cover $X' \to X$, the map \[ \kappa_{X'/k} \colon X'(k) \to \seck{X'} \] is injective. The injectivity of the extended non-abelian Kummer map \[ \kappa_{X'/k} \colon X'(k) \amalg \cC_{X'} \to \seck{X'} \] follows from the injectivity of the extended non-abelian Kummer map on $X$.
\end{rem}
\subsection{Geometric Fibrations}

\begin{defn}[{\cite[Definition 11.4]{fr82}}]\label{defn:geom_fib}
A \emph{special geometric fibration} is a morphism $f \colon X \to S$ of schemes fitting into a commutative diagram:
$$
\xymatrix{
X  \ar@{^{(}->}[r]^{j} \ar[rd]^f & \overline{X} \ar[d]^{\overline{f}} & Y\ar@{_{(}->}[l]^{i} \ar[ld]^g \\
\empty & S &\empty
}
$$
satisfying the following conditions:
\begin{enumerate}
\item $i$ is a closed embedding
\item $j$ is an open immersion which is dense in every fiber of $\overline{f}$, and $X = \overline{X} - Y$
\item $\bar{f}$ is smooth and proper
\item\label{item:intersection_smooth}  $Y$ is a union of schemes $Y_1,\cdots,Y_m$, with $Y_i$ of pure relative codimension $c_i$ in $\cl{X}$ over $S$, with the property that every intersection $Y_{i_1} \cap \cdots \cap Y_{i_k}$ is smooth over $S$ of pure codimension $c_{i_1} + \cdots + c_{i_k}$.
\end{enumerate}

More generally, a \emph{geometric fibration} is a map $f \colon X \to S$ admitting a Zariski covering $\{V_i \to S\}$ such that $f_{V_i} \colon X \times_S V_i \to V_i$ is a special geometric fibration for all $i$.

\end{defn}

\begin{defn}\label{defn:cusp_fib}
Let $k$ be a field of characteristic zero and $S$ a $k$-variety. Let $f \colon X \to S$ be a geometric fibration with geometrically connected fibers. Suppose we are given cuspidal data $(S,\cC_S)$ for $S$ and $(X_s,\cC_{X_s})$ for the fiber $X_s$ above every $s \in S(k)$. 
For $s \in S(k)$, let $\iota_s \colon X_s \to X$ denote the inclusion. We define $\cC_{X,f}$ by
$$\cC_{X,f} = f_*^{-1}(\cC_S) \cup \left(\bigcup_{s \in S(k)} (\iota_s)_*(\cC_{X_s})\right) \subseteq \seck{X}.$$ We call $(X,\cC_{X,f})$ the \emph{cuspidal datum induced by $f$}.

\end{defn}
Recall that a $k$-variety $X$ is an étale $K(\pi,1)$ (c.f. \cite[\S 2.3]{schsx16}) if the profinite \'etale homotopy groups $\pi_i(X,\overline{x})$ vanish for $i \geq 2$ and some (equivalently, any) geometric point $\overline{x}$.

\begin{prop}\label{prop:from_onpoonence}
Let $f \colon X \to S$, $\cC_S$, $\cC_{X_s}$, and $\cC_{X,f}$ be as in Definition \ref{defn:cusp_fib}. Assume that $S$ is normal, that $S$ is a $K(\pi,1)$, and that all geometric fibers of $f$ are étale $K(\pi,1)$'s. Then $X$ is a $K(\pi,1)$ as well. In that case, if $(S,\cC_S)$ and all $(X_a,\cC_{X_a})$ for $a \in S(k)$ satisfy injectivity in the section conjecture, then so does $(X,\cC_{X,f})$. The same is true for the surjectivity in the section conjecture.
\end{prop}
\begin{proof}
The first part follows from the long exact sequence of \cite[Theorem 11.5]{fr82} for $L$ the set of all primes; see also \cite[Proposition 2.8]{schsx16}. 

Assume now that $(S,\cC_S)$ and $(X_a,\cC_{X_a})$, where $a \in S(k)$, satisfy injectivity in the section conjecture. Then to show that $(X,\cC_{X,f})$ satisfies injectivity in the section conjecture, we use Remark \ref{remark:cuspidalequiv}. By Remark \ref{remark:cuspidalequiv}, $Z_{\pi_1(S^s)}(\kappa_{S/k}(a)) = 1$ for every $a \in S(k)$. The long exact sequence of \cite[Proposition 31]{sxbook13} shows that $(\iota_a)_* \colon \seck{X_a} \to \seck{X}$ is injective. This injectivity, together with the injectivity of the non-abelian Kummer map on $S$ and $X_a$ for $a \in S(k)$ shows that $X(k) \to \seck{X}$ is injective. The triviality of the centralizer $Z_{\pi_1(X^s)}(\kappa_{X/k}(b))$ for $b \in X(k)$ follows from the exact sequence of centralizers supplied by \cite[Proposition 31]{sxbook13}. 

Lastly, to show that the extended Kummer map is injective, it suffices to prove that $\kappa_{X/k}(X(k)) \cap \cC_{X,f} = \emptyset$. If $\kappa_{X/k}(b) \in \cC_{X,f}$ where $b \in X(k)$, we find by the injectivity of the extended Kummer map on $S$ that $\kappa_{X/k}(b) \in (\iota_s)_*(\cC_{X_s})$ for $s = f(b)$. However, from the injectivity of $(\iota_s)_* \colon \seck{X_s} \to  \seck{X}$ it follows that $\kappa_{X_s/k}(b) \in \cC_{X_s}$, which contradicts the injectivity of the extended Kummer map on $X_s$. 

Assume now that $(S,\cC_S)$ and $(X_a,\cC_{X_a})$ for every $a \in S(k)$ satisfy surjectivity in the section conjecture. Let $s \in \seck{X}$ be a section. If $f_*(s) \in \cC_S$, it is clear that $s \in \cC_{X,f}$. Otherwise, surjectivity of the section conjecture on $S$ shows that $f_*(s) = \kappa_{S/k}(b)$ for some $b \in S(k)$. Then \cite[Proposition 31]{sxbook13}, together with the paragraph preceding it, shows that $s \in (\iota_b)_*(\seck{X_b})$. The surjectivity of the section conjecture on $X_b$ allows us to conclude.
\end{proof}

\subsection{Cuspidal data of hyperbolic polycurves and the section conjecture}
Let $k$ be a field of characteristic zero. A \emph{hyperbolic polycurve} over $k$ is a $k$-variety $X$ which admits a sequence of geometric fibrations 
\[ X = X_n \to X_{n-1} \to \cdots \to X_1 \to X_0 = \Spec{k} \]
of $k$-varieties such that the geometric fibers of each $X_i \to X_{i-1}$ are all hyperbolic curves. We call such a geometric fibration $X_i \to X_{i-1}$ a \emph{relative hyperbolic curve} (over $X_{i-1}$). Such a sequence for $X$ is called a \emph{sequence of parameterizing morphisms}. Then we can inductively define a cuspidal datum $C_{X,f}$ of $X= X_n$ via Definition \ref{defn:cusp_fib}, using the standard cuspidal datum on the fibers.
\begin{rem}
Our definition of hyperbolic polycurves is equivalent to that of \cite[Definition 2.1(ii)]{Hoshi14}. It suffices to show that $X_i \to X_{i-1}$ is a relative hyperbolic curve in the sense of \cite[Definition 2.1(i)]{Hoshi14}. Since the base $X_{i-1}$ is normal, compactifications are unique by \cite[Definition 2.2]{Hoshi14}, so by gluing compactifications, we may assume that $X_i \to X_{i-1}$ is a special geometric fibration. Then since $X_i \to X_{i-1}$ is of relative dimension one, the boundary is finite étale over $X_{i-1}$ by (\ref{item:intersection_smooth}) of Definition \ref{defn:geom_fib}, and it follows immediately that $X_i \to X_{i-1}$ is a relative hyperbolic curve in the sense of Hoshi. That any hyperbolic polycurve in the sense of Hoshi is a hyperbolic polycurve in our sense is immediate. 
\end{rem}
The following proposition studies the independence of $\cC_{X,f}$ from the given set of parameterizing morphisms.
\begin{prop} \label{prop:equivsection}
Let $k$ be a field of characteristic zero, and assume that every hyperbolic curve satisfies injectivity in the section conjecture for the standard cuspidal datum. Then the cuspidal datum $(X,\cC_{X,f})$ is independent of the sequence of parameterizing morphisms for every hyperbolic polycurve $X$ over $k$ if and only if the section conjecture holds for hyperbolic curves over $k$.
\end{prop}
\begin{proof}
Assume first that the section conjecture holds over $k$. Given a hyperbolic polycurve $X$, we note that by Proposition \ref{prop:from_onpoonence}, $(X,\cC_{X,f})$ satisfies injectivity as well as surjectivity in the section conjecture. Thus, $\cC_{X,f} = \seck{X} \setminus \kappa_{X/k}(X(k))$; this set is clearly independent of the sequence of parameterizing morphisms.

Conversely, assume the independence of the cuspidal datum for all hyperbolic polycurves over $k$. Let $S$ be a hyperbolic curve over $k$ and set $Y = \mathbb{P}^1_k-\{0,1,\infty\}$. We let $X = S \times_k Y$ and observe that $X$ is a hyperbolic polycurve via both the first and second projection. To simplify notation we set $f = \pr_1$ and $g = \pr_2$. Give both $S$ and $Y$ the standard cuspidal data $(S,C_S)$ and $(Y,C_Y)$ as in Definition \ref{defn:cusp_curve}. Note that $\seck{X} = \seck{S} \times \seck{Y}$, and that under this identification we have
\[ \cC_{X,f} = C_S \times \seck{Y} \cup \left( \kappa_{S/k}(S(k)) \times C_Y  \right) \subset \seck{X} \]  and 
\[\cC_{X,g} = \seck{S} \times C_Y \cup \left( C_S \times \kappa_{Y/k}(Y(k)) \right) \subset \seck{X}.  \] If $\cC_{X,f} = \cC_{X,g}$, then for an arbitrary $s \in \seck{S}$ and $b \in C_Y$, the product $s \times b \in \cC_{X,g}$. By the definition of $\cC_{X,f}$ it follows that $s$ is either cuspidal or rational, i.e., the surjectivity of the section conjecture holds for $S$. Since injectivity further holds, the section conjecture holds for $S$.
\end{proof}
For a hyperbolic polycurve $X$, we relate the (closed) cuspidal set $C_X$ constructed in Theorem \ref{thm:cuspidals_independent_of_compactification} to the closure of $C_{X,f}$.

\begin{prop}\label{prop:fibinclusion}
Let $k$ be a field of characteristic zero, and let $X$ be a hyperbolic polycurve over $k$. Then $C_{X} \subset \overline{\cC_{X,f}}$.
\end{prop}

\begin{proof}

We prove by induction on the length of the tower of parameterizing morphisms that there exists an snc-compactification $Z$ of $X$ such that $C_{X,Z} \subset \cC_{X,f}$. By taking closures, we obtain our Proposition. Obviously the desired result holds when $X$ is a hyperbolic curve.

Let $X$ be a hyperbolic polycurve of dimension greater than one, and let $f \colon X \to S$ be the first map in the sequence of parameterizing morphisms, and write $\cC_{f'}$ for the cuspidal set associated to the truncated hyperbolic polycurve structure on $S$. By induction, choose an snc-compactification $S \hookrightarrow T$ such that $C_{S,T} \subset \cC_{f'}$. 

We let $\cl{f} \colon \cl{X} \to S$ denote the smooth compactification of $f$ with boundary $Y$ (c.f. Definition \ref{defn:geom_fib}). Nagata's compactification theorem gives a factorization $\cl{X} \to Z \to T$, where $\cl{X} \to Z$ is a schematically dense open immersion and $Z \to T$ is proper. Note that $\overline{X} \to Z \times_T S$ is an isomorphism since $\cl{f}$ is proper. By \cite[Theorem 36]{KollarSeattle}, we can find a proper birational map $r \colon Z' \to Z$ with $Z'$ smooth and with the property that $r$ is an isomorphism over $\cl{X}$. By blowing up we can assume that $Z' \setminus X$ is a simple normal crossings divisor \cite[Theorem 35]{KollarSeattle}. Thus, by replacing $Z$ by $Z'$ we can assume that there is a smooth and proper compactification 
\[ j \colon X \to Z \] of $X$ such that $Z \setminus X$ is a simple normal crossings divisor, and an extension $p \colon Z \to T$ of $X \to T$. Further, since $Y$ is smooth, the further blow-ups can be chosen so that they are isomorphisms over $\cl{X}$, by the smooth functoriality properties of \cite[Theorem 35]{KollarSeattle}. This shows that $Z \times_T S = \cl{X}$, and that $Z \times_T S \to S$ is smooth with relative snc boundary.

To prove that $C_{X,Z} \subset \cC_{X,f}$, let $y \in (Z \setminus X)(k)$. Let $X_{(y)}$ be the scheme of nearby points and let $s \in \seck{X_{(y)}}$. Let $i \colon X_{(y)}\to X$ be the natural map and denote by $i_*(s) \in \seck{X}$ the natural section. Let $b \colonequals p(y) \in T(k)$. If $b \in (T \setminus S)(k)$, then by functoriality 
\[ f_*(i_*(s)) \in C_{S,T} \subset \cC_{f'}. \] Thus, in this case, $i_*(s) \in \cC_{X,f}$. 

On the other hand, if $b \in S(k)$, note that $y$ is a boundary point of the smooth compactification $\cl{X}_b$ of $X_b$. We let $(X_b)_{(y)}$ be the scheme of nearby points of $X_b$ at $y$. There is a natural map $(X_b)_{(y)} \to X_{(y)}$, and we claim that the induced map 
\[ \pi_1((X_b)_{(y)}) \to \pi_1(X_{(y)}) \] is an isomorphism over $G_k$. This follows from Proposition \ref{prop:fundamentalofU}: since $Y \to S$ is finite étale, the local components that pass through $y$ restrict bijectively and with multiplicity one to those of $X_b$ passing through $y$.  
This shows that $i_*(s)$ lies in the image of $C_{X_b}$ in $\seck{X}$, thus $i_*(s) \in \cC_{X,f}$.  Thus, $C_{X,Z} \subset \cC_{X,f}$.
\end{proof}

{\def\section*#1{}

\bibliography{references,etale_homotopy,chabauty-kim,ag,logcusp}
}

\end{document}